\documentclass[12pt,draft,amsmath,amscd,amsbsy,amssymb]{amsart}
 \relax

\chardef\bslash=`\\ 

\makeatletter \def\verbatim{\interlinepenalty\@M \@verbatim
\leftskip\@totalleftmargin\advance\leftskip2pc
\frenchspacing\@vobeyspaces \@xverbatim} \makeatother \hfuzz1pc

\makeatletter \def\dgt@k{\dg@DX=-3 \dg@DY=2 \dg@SIZE=3}
\makeatother

\makeatletter \def\dgt@kk{\dg@DX=3 \dg@DY=-1 \dg@SIZE=3}

\theoremstyle{plain} \newtheorem{thm}{Theorem}[section]
\newtheorem{cor}[thm]{Corollary}

\theoremstyle{definition} 
 \newtheorem{que}[thm]{Question}

\newcommand{\Real}{\mathbb R}

\newcommand{\cw}{\operatorname{\mathrm{cw}}}
\newcommand{\pcw}{\operatorname{\mathrm{pcw}}}
\newcommand{\crw}{\operatorname{\mathrm{crw}}}
\newcommand{\cc}{\operatorname{\mathrm{cc}}}

\newcommand{\pr}{\mathrm{pr}}

\usepackage[all]{xy}

\begin{document}

\title[On convex bodies of constant width]
{On convex bodies of constant width}
\author{L.E. Bazylevych}
\address{National University ``Lviv Polytechnica", 12 Bandery Str.,
79013 Lviv, Ukraine}

\author{M.M. Zarichnyi}
\address{Lviv National University, 1 Universytetska Str.,
79000 Lviv, Ukraine}
\email{mzar@litech.lviv.ua}
\thanks{}
\subjclass{54B20, 52A20, 46A55}
\keywords{Convex body, constant width, $Q$-manifold, soft map.}
\dedicatory{Dedicated to Professor E.D. Tymchatyn on the occasion of
his 60th anniversary}


\begin{abstract} We present an alternative proof of the following
fact: the hyperspace of compact closed subsets of constant width in
$\Real^n$ is a contractible Hilbert cube manifold. The proof also
works for certain subspaces of compact convex sets of constant
width as well as for the pairs of compact convex sets of constant
relative width. Besides, it is proved that the projection map of
compact closed subsets of constant width is not 0-soft in the sense
of Shchepin, in particular, is not open.
\end{abstract}

\maketitle

The topology of the hyperspace of compact convex sets in euclidean
spaces was investigated by different authors; see, e.g.
\cite{NQS,M,B,B1}.

In this note we consider some topological properties of (the maps
of) compact convex bodies of constant width. A convex set in
euclidean space is said to be of {\em constant width\/} $d$ if the
distance between two supporting hyperplanes equals $d$ in every
direction. To be more formal, denote by $h_K\colon S^{n-1}\to
\Real$ the {\em support function} of a convex body $K$ in $\Real^n$
defined as follows: $h_K(u)=\max\{\langle x,u\rangle\mid x\in K\}$.
Here, as usual, $\langle ,\rangle$ stands for the standard inner
product in $\Real^n$ and $S^{n-1}$ is the unit sphere in $\Real^n$.
The {\em widths function} of $K$ is the function $w_K\colon S^{n-1}\to
\Real$ defined by the formula $w_K(u)=h_K(u)-h_K(-u)$. A convex
body $K$ is of constant width $d>0$ provided $w_K$ is a constant
function taking the value $d$.

It was proved by the first-named authors \cite{B} that the
hyperspace of compact convex bodies of constant (non-specified)
width in euclidean space of dimension $\ge2$ is homeomorphic to the
punctured Hilbert cube. In Section 1 we present a more direct proof
of this result. The technique allows also to prove that some
subspaces of the above mentioned hyperspace as well as the
hyperspace of pairs of compact convex sets of constant relative
width are manifolds modeled on the Hilbert cube ($Q$-manifolds).

In connection to the topology of the hyperspace of compact convex
sets of nonmetrizable compact subsets in locally convex spaces it
was proved in \cite{IZ} that, for any affine continuous onto map of
convex subsets in metrizable locally convex spaces, the natural map
of the hyperspaces of compact convex subsets is soft in the sense
of Shchepin \cite{S}. In Section 2 we demonstrate that this is not
the case if we restrict ourselves with the compact convex sets of
constant width; the map under consideration is not even open.

\section{Hyperspaces of compact convex bodies of constant width}

By $\cc(\Real^n)$ we denote the hyperspace of compact convex
subsets in $\Real^n$. We equip $\cc(\Real^n)$ with the topology
generated by the Hausdorff metric.

By $\cw_d(\Real^n)$ we denote the subset of $\cc(\Real^n)$
consisting of convex bodies of constant width $d>0$ in $\Real^n$.
We also put $\cw_0(\Real^n)=\{\{x\}\mid x\in \Real^n\}$. In
\cite{B} it is proved that the hyperspace
$\cw(\Real^n)=\cup_{d>0}\cw_d(\Real^n)$, $n\ge2$, is a contractible
Hilbert cube manifold. Here we essentially simplify the proof of
this result.

\begin{thm}\label{t:1} Let $D\subset[0,\infty)$ be a convex subset such that
$D\cap(0,\infty)\neq\emptyset$. The hyperspace
$\cup\{\cw_d(\Real^n)\mid d\in D\}$ is a contractible Hilbert cube
manifold.
\end{thm}

\begin{proof}

First, we embed $\cw(\Real^n)$ as a convex subset of a Banach
space.

Define a map $\varphi\colon\cw(\Real^n)\to C(S^{n-1})$ by the
formula $\varphi(K)(x)=h_K(x)$, $K\in \cw(\Real^n)$. It is a
well-known fact that $\varphi$ is a continuous map. Moreover, it is
obvious that $\varphi$ is an embedding which is an affine map in
the sense that $\varphi(tA+(1-t)B)=t\varphi(A)+(1-t)\varphi(B)$ for
every $A,B\in
\cw(\Real^n)$ and $t\in[0,1]$. The image of $\varphi$ is a locally
compact convex subset of $C(S^{n-1})$.

We are going to prove that, for any $d>0$, the space $\cw(\Real^n)$
is infinite-dimensional.

First consider the case $n=2$. Let $K$ denote the Reuleaux triangle
in $\Real^2$ that is the intersection of the closed balls of radius
$d$ centered at $(0,0)$, $(d,0)$, and $(d/2, d\sqrt{3}/2)$. For any
$\alpha\in[0,2\pi]$, denote by $K_\alpha$ the convex body obtained
by rotation of $K$ by angle $\alpha$ conterclockwise around the
origin. We show that the the set
$\{K_\alpha\mid\alpha\in[0,2\pi)\}$ contains linearly independent
subset of arbitrary finite cardinality. To this end, one has to
demonstrate that the family of the support functions
$h_{K_\alpha}\colon S^1\to
\Real$ contains linearly independent subset of arbitrary finite
cardinality.

We identify $S^1$ with the set $\{e^{it}\mid t\in[0,2\pi]\}$. It is
easy to see that
$h_{K_\alpha}(e^{i\gamma})=h_{K}(e^{i(\gamma-\alpha)})$. Elementary
geometric arguments demonstrate that
$$h_K|[0,\pi/3]=1,\ h_K|[\pi,4\pi/3]=0,\
h_K((\pi/3,\pi)\cup(4\pi/3,2\pi))\subset (0,1).$$

Fix a natural number $k$. For each $j=0,1,\dots,k-1$, let
$h_j=h_{K_{(j\pi)/(3k)}}$. In order to demonstrate that the
functions $h_j$, $j=0,1,\dots,k$, are linearly independent,
consider a linear combination $g=\sum_{j=0}^{k}\lambda_jh_j$.
Suppose that $g=0$, then
\begin{align*}
   & g(\pi/3)=\sum_{j=0}^{k}\lambda_jh_j(\pi/3)=
\sum_{j=0}^{k}\lambda_j=0, \\
   & g((\pi/3)+(\pi/(3k)))=
\lambda_0h_0((\pi/3)+(\pi/(3k)))\\ +&\sum_{j=1}^{k}\lambda_j
h_j((\pi/3)+(\pi/(3k)))=0,
\end{align*}
whence $\lambda_0=0$. Consequently evaluating the function $g$ at
the points $(\pi/3)+(j\pi/(3k))$, $j=2,\dots, k$, we conclude that
$\lambda_i=0$ for every $i=0,1,\dots,k$.

In the case $n>2$, consider the family of $n$-dimensional simplices
in $\Real^n$ all whose edges are of length $d$. For every such a
simplex, $\Delta$, consider the intersection of all closed balls of
radius $d$ centered at the vertices of $\Delta$ and containing
$\Delta$. The obtained set, $L$, is obviously a compact convex
subset in $\Real^n$ of constant width $d$. Denote by $\mathcal L$
the family of all compact convex subsets in $\Real^n$ that can be
obtained in this way. Let $\pr\colon \Real^n\to\Real^2$ denote the
projection. This projection generates the map
$\cc(\pr)\colon\cc(\Real^n)\to\cc\Real^m$, $f(A)=\pr(A)$. The map
$\cc(\pr)$ is an affine map.

Then $\{\cc(\pr)(L')\mid L'\in\mathcal L\}\supset
\{K_\alpha\mid\alpha\in[0,2\pi]\}$ and we conclude that the space
$\mathcal L$ and, therefore $\cw_d(\Real^n)$ is
infinite-dimensional.

To finish the proof, apply the results on topology of metrizable
locally compact convex subsets in locally convex spaces \cite{BP}.
Since the space $\cw(\Real^n)$ is easily shown to be locally
compact, the Keller theorem (see \cite{BP}) implies that
$\cw(\Real^n)$ is a $Q$-manifold.

\end{proof}

\begin{cor} Let $d_0\ge0$. The hyperspace $\cup\{\cw_d(\Real^n)\mid d\ge d_0\}$
is homeomorphic to a punctured Hilbert cube $Q\setminus\{*\}$.
\end{cor}
\begin{proof}
By a result of Chapman \cite{C}, it is sufficient to prove that
there exists a proper homotopy $$H\colon \cup\{\cw_d(\Real^n)\mid
d\ge d_0\}\times[1,\infty)\to \cup\{\cw_d(\Real^n)\mid d\ge
d_0\}.$$ This homotopy can be defined in an obvious way,
$H(K,t)=tK$, $(K,t)\in \cup\{\cw_d(\Real^n)\mid d\ge
d_0\}\times[1,\infty)$.
\end{proof}

\begin{thm}\label{t:convex}
Let $X$ be a convex subset in $\Real^n$ that contains a closed
square of side $d$. Then for any convex subset $D$ of $[0,d]$ with
$D\cap(0,d]\neq\emptyset$ the set
$\cw_D(X)=\cc(X)\cap\cw_D(\Real^n)$ is a Hilbert cube manifold.
\end{thm}

\begin{proof}
Note that the family $\{L-\{x\}\mid L\in\cw_D(X),\ x\in X\}$
contains the family $\{K_\alpha\mid\alpha\in[0,2\pi)\}$ defined in
the proof of Theorem \ref{t:1} (for fixed $d$). Since $X$ is a
subset of a finite-dimensional linear space and
$\{K_\alpha\mid\alpha\in[0,2\pi)\}$ contains an infinite linearly
independent family, we conclude that $\cw_D(X)$ also contains an
infinite linearly independent family. Therefore, the space
$\cw_D(X)$ is infinite-dimensional. Then we apply the arguments of
the proof of Theorem \ref{t:1}
\end{proof}

Remark that it follows, in particular, from Theorem \ref{t:convex}
that the hyperspace of convex bodies that can be rotated inside a
square is homeomorphic to the Hilbert cube, because this hyperspace
is a compact contractible $Q$-manifold (see \cite{Ch}).

\section{On softness of the projection map}

A map $f\colon X\to Y$ is {\it soft} (respectively {\em $n$-soft})
if for every commutative diagram
$$\xymatrix{A \ar[r]^{\psi}\ar[d]_i & X\ar[d]^{f} \\
Z\ar[r]_{\varphi} & Y,}$$ where $i\colon A\to Z$ is a closed
embedding into a paracompact space $Z$ (respectively a paracompact
space $Z$ of covering dimension $\le n$), there exists a map
$\Phi\colon Z\to X$ such that $\Phi | A=\psi$ and $f\Phi=\varphi$.
The notion of ($n$-)soft map was introduced by E.V. Shchepin
\cite{S}.

Let $\pr\colon \Real^n\to\Real^m$ denote the projection, $n\ge m$.
As we already remarked, this projection generates the map
$\cc(\pr)\colon\cc(\Real^n)\to\cc\Real^m$ and by
$$p=\cw(\pr)\colon\cw(\Real^n)\to\cw\Real^m$$ we denote its
restriction onto $\cc(\pr)|\cw(\Real^n)$. It is proved in \cite{IZ}
that the map $\cc(\pr)$ is soft.

The images of the fibers of the map $p$ under the embedding
$\varphi$ are obviously convex. It is natural to ask whether the
map $p$ is also soft. As the following result shows, the answer
turns out to be negative. The idea of the proof is suggested by S.
Ivanov.
\begin{thm}\label{t:soft}
The map $p\colon\cw(\Real^3)\to\cw\Real^2$ is not 0-soft.
\end{thm}
\begin{proof}
Consider the compactum $L$ in $\Real^3$ which is the intersection
of the closed balls of radius 2 centered at the points $(0,1,0)$,
$(0,-1,0)$, $(-1,0,\sqrt{2}$, and $(1,0,\sqrt{2}$ (a Reuleaux
tetrahedron). It is well-known that $L$ is of constant width and
$f(L)=K$, where $K$ denotes the disc of radius 1 centered at the
origin in $\Real^2$. For every $i$, denote by $K_i$ the compactum
in $\Real^2$ described as follows. Let
$$x_j=(\cos(2\pi j/(2i+1)),\sin(2\pi j/(2i+1))),\ j=0,1,\dots,2i.$$
The compactum $K_i$ is the intersection of the discs of radius
$\|x_0-x_i\|$ centered at the points $x_j$, $j=0,1,\dots,2i$.
Obviously, $K_i\in
\cw\Real^2$ and $\lim_{i\to\infty}K_i= K$ in the Hausdorff metric.

Let $S=\{0\}\cup\{1/n\mid n\in\mathbb N\}$, $f\colon\{0\}\to
\cw(\Real^3)$ be the map that sends 0 into $L$, and
$F\colon S\to\cw(\Real^2)$ be the map such that $F(0)=K$ and
$F(1/n)=K_n$, $n\in\mathbb N$.

Suppose now that the map $p$ is 0-soft. Then there exist a map
$G\colon S\to \cw(\Real^3)$ such that $G(0)=f(o)=L$ and $pG=F$.
Denoting $G(1/i)$ by $L_i$, we see that $p(L_i)=K_i$ and
$\lim_{i\to\infty}L_i=L$. Since every diameter of a convex body
(i.e., a segment that connects the points at which two parallel
supporting hyperplanes touch the body) of constant width is
orthogonal to the supporting planes, for every $i$, there exists a
plane in $\Real^3$ containing the set $\pr^{-1}(\partial K_i)\cap
L_i$. From this we conclude that there exists a plane in $\Real^3$
containing the set
$$\pr^{-1}(\partial K)\cap L\supset \{(0,1,0), (0,-1,0),
(-1,0,\sqrt{2},(1,0,\sqrt{2}\},$$ a contradiction.
\end{proof}
\begin{cor}
The map $f\colon\cw(\Real^3)\to\cw\Real^2$ is not open.
\end{cor}

\section{Pairs of compact convex bodies of constant relative width}

Two convex bodies $K_1,K_2\subset
\Real^n$ are said to be a {\em pair of constant width} if
$K_1 -K_2$ is a ball (Maehara \cite{Ma}). We denote by
$\crw(\Real^n)$ the set of all pairs of compact convex bodies of
constant relative width. The set $\crw(\Real^n)$ is topologized
with the subspace topology of $\cc(\Real^n)\times\cc(\Real^n)$.
Note that the space $\cw(\Real^n)$ can be naturally embedded into
$\crw(\Real^n)$ by means of the mapping $A\mapsto (A,A)$.

\begin{thm}
The space $\crw(\Real^n)$ is a contractible $Q$-manifold.
\end{thm}

\begin{proof}
First, note that the map $\beta\colon (A,B)\mapsto (h_A,h_B)$
affinely embeds the space $\crw(\Real^n)$ into $C(S^{n-1})\times
C(S^{n-1})$ and the image of this embedding is a convex subset of
$C(S^{n-1})\times C(S^{n-1})$. Since the space $\cw(\Real^n)$ is
infinite-dimensional, so is $\crw(\Real^n)$. It is easy to see that
$\cw(\Real^n)$ is locally compact. Arguing like in the proof of
Theorem \ref{t:1}, we conclude that $\cw(\Real^n)$ is a
contractible $Q$-manifold.
\end{proof}

\section{Remarks and open questions}

The set $\cw_d(\Real^n)$ is the preimage of the set of closed balls
of radius $d$ under the {\em central symmetry map}
$c\colon \cc(\Real^n)\to
\cc(\Real^n)$, $c(A)=(A-A)/2$. Since $\cw_d(\Real^n)$ is an absolute retract, this
directly leads to the following question.
\begin{que}
Is the central symmetry map soft as a map onto its image.
\end{que}

Let $L$ be a Minkowski space (i.e. a finite-dimensional Banach
space) of dimension $\ge2$. Given a compact convex body $K\subset
L$, one says that $K$ is of constant width if the set $K-K$ is a
closed ball centered at the origin. A natural question arizes: for
which spaces $L$ the results of this paper can be extended over the
hyperspaces of convex bodies of constant width in $L$? Note that
for the space $(\Real^n,\|\cdot\|_\infty)$ the hyperspace of convex
bodies of constant width coincides with that of closed balls in it
and therefore is finite-dimensional. In addition, a counterpart of
Theorem \ref{t:soft} does not hold for this space.

\begin{que}
Describe the topology of the hyperspace of smooth convex bodies of
constant width.
\end{que}

Every compact convex body of constant width is a rotor in cubes as
well as in some other polyhedra.

\begin{que}
Find counterparts of the results of this paper for rotors in
another polyhedra (e.g., equilateral triangles).
\end{que}


\end{document}